\documentclass[12pt,a4paper]{article}
\usepackage{amsmath}
\usepackage{amsthm}
\usepackage{amssymb}
\usepackage{appendix}
\usepackage{xcolor}
\usepackage{hyperref}

\numberwithin{equation}{section}
\numberwithin{figure}{section}

\theoremstyle{plain}
\newtheorem{thm}{Theorem}[section]
\newtheorem{prop}[thm]{Proposition}
\newtheorem{lem}[thm]{Lemma}

\theoremstyle{remark}

\newtheorem*{rem*}{Remark}

\theoremstyle{definition}
\newtheorem{defn}[thm]{Definition}

\usepackage{setspace}
\global\long\def\apa{\alpha}%
 
\global\long\def\ba{\beta}%
 
\global\long\def\ga{\gamma}%

\global\long\def\da{\delta}%
 
\global\long\def\Da{\Delta}%

\global\long\def\oa{\omega}%
 
\global\long\def\Oa{\Omega}%
 
\global\long\def\vn{\varepsilon}%
 
\global\long\def\pl{\partial}%
 
\global\long\def\na{\nabla}%

\global\long\def\fc#1#2{\frac{#1}{#2}}%

\global\long\def\ol#1{\overline{#1}}%

\global\long\def\wt#1{\widetilde{#1}}%
 
\global\long\def\mb#1{\mathbb{#1}}%
 
\global\long\def\ml#1{\mathcal{#1}}%

\begin{document}

\title{The enclosure method for semilinear elliptic equations with power-type nonlinearities}
\author{Rulin Kuan \\
National Cheng Kung University, Tainan, Taiwan\\
Email: rkuan@ncku.edu.tw}
\date{}
\maketitle

\begin{abstract}
We employ the enclosure method to reconstruct unknown inclusions within an object
that is governed by a semilinear elliptic equation with power-type nonlinearity. 
Motivated by \cite{LLLS2021}, we tried to solve the problem without using special solutions
such as complex geometrical optics solutions to the equation. Instead, we construct 
approximate solutions obtained from Taylor approximation of the solution operator.
By incorporating such approximate solutions into the definition of an indicator functional, 
we are able to use the classical Calder\'{o}n-type harmonic functions to 
accomplish the reconstruction task.\\

\end{abstract}

\section{Introduction}

The enclosure method is a method to reconstruct the shapes of unknown inclusions within an object 
from boundary measurements. It is initiated by Ikehata \cite{I2000, I1999-2}, and has been successfully 
applied to a variety of equations (see e.g. \cite{UW2008, zhou2009reconstructing, 
nagayasuUW2011, K2012, KLS2015, I2022, I2023}). However, most studies so far consider 
linear equations, and few consider nonlinear ones (see however \cite{BKS2015} and \cite{BHKS2018}, which consider the $p$-Laplace equation.)
Therefore, it would be of interest to investigate the applicability of the enclosure method to a wider range of nonlinear equations.

\subsection{The problem}\label{sec:prob}

In this paper, we consider the following semilinear elliptic equation with power-type nonlinearity:
\begin{align}\label{eq:main}
\Delta u+q(x)u^{m}=0\quad\mbox{in}\quad\Oa,
\end{align}
where $q\in L^\infty(\Oa;\mb R)$; $m\in\mb N$, $m\ge 2$; 
and $\Oa$ is a bounded smooth domain in $\mb R^n$, $n\ge 2$.
We assume 
\begin{align*}
q(x)=q_0(x)+\chi_{D}(x)q_{D}(x),
\end{align*}
where $D\subset \Oa$ is an unknown open set (not necessarily connected), 
and $q_0, q_D\in L^\infty(\Oa;\mb R)$.
We assume that $q_0$ is known, which represents the 
expected background coefficient. On the other hand,  
$q_D$ (as well as the position of $D$) accounts for the presence of some 
``obstacles'' and is unknown. Nevertheless, we assume $q_D$ satisfies the following jump condition
so that the obstacles can be distinguished from the background material:
\begin{align}\label{eq:jump}
\begin{aligned}
&q_D(x)\ge \mu\mbox{ for all } x\in D\mbox{ or }q_D(x)\le -\mu\mbox{ for all }x\in D,\\
&\mbox{where }\mu>0 \mbox{ is a constant}.
\end{aligned}
\end{align}
Our goal is to apply the enclosure method to reconstruct the convex hull of $D$.
(See Theorem \ref{thm:main_thm} for our main result.)

\subsection{An exemplified result of the enclosure method}\label{sec:in2en}

Before diving into Equation \eqref{eq:main}, let us give a classical result
(restated from Theorem 3.1 in \cite{I1999-2})
to illustrate the key elements of the enclosure method. Here we ignore all the regularity 
problems, and just assume that all functions and domains are ``regular enough'' for convenience.
Consider the following conductivity equation with homogeneous background:
\begin{align}\label{simple_ex}
\na\cdot(1+\chi_D \ga_D(x))u=0\quad\mbox{in}\quad\Oa, 
\end{align}
where $D\subset \Oa\subset\mb R^n$, $n\ge 2$, and $\ga_D:\Oa\to\mb R$ 
satisfies $\ga_D(x)>\mu$ for some positive constant $\mu$. 

For a Dirichlet boundary condition ``$u=f$ on $\pl\Oa$'' assigned to \eqref{simple_ex}, 
let $u_f$ denote the corresponding solution.
Similarly, let $u_{0,f}$ denote the solution for the case $\ga_D=0$ (or equivalently, $D=\emptyset$).
For a function $u$ defined in $\Oa$, let $\pl_\nu u$ be the outward normal derivative of $u$ on $\pl\Oa$.
Then,
define the following ``indicator functional'' on Dirichlet data $f$:
\begin{align*}
E(f)=\int_{\pl\Oa}(\pl_\nu u_f -\pl_\nu u_{0,f}) \ol{f}\,d\sigma, 
\end{align*}
and the following class of ``test (Dirichlet) data'':
\begin{align*}
f_h(x)=e^{-\frac{1}{h}(x\cdot \oa+ix\cdot\oa^{\perp}-t)}|_{\pl\Oa},
\end{align*}
where $h>0$ is a small parameter, $t\in\mb R$, and $\oa,\oa^\perp\in\mb{S}^{n-1}$ (the unit sphere in $\mb R^n$ centered at 
the origin) are such that $\oa\perp\oa^\perp$.


Then it can be proved that $\lim_{h\to 0}E(f_h)=0$ if and only if $t<\inf_{x\in D}x\cdot \oa$, or 
equivalently, the half-space $\{x\cdot\oa \le t\}$ does not intersect $\ol{D}$. 
Thus, by varying the direction $\oa$ and the parameter $t$, and then examining the limiting 
behavior of $E(f_h)$, we can theoretically ``enclose'' the region $D$ by half-spaces, 
and therefore reconstruct the convex hull of $D$. 

Note that functions of the form $\apa e^{\beta x\cdot(\oa+i\oa^\perp)}$, 
which provide the test data in the above reconstruction, are harmonic functions. 
In other words, they are solutions to the background equation for \eqref{simple_ex}. 
Such complex harmonic functions were employed in Calder\'{o}n's pioneering work \cite{Cal2006},
and are prototypes of the so-called complex geometrical optics (CGO) solutions, introduced in \cite{SU1987}.

\subsection{Motivation and innovation}\label{sec:monin}

As illustrated by the example above, to apply the enclosure method one needs a suitable 
indicator functional and a suitable class of test data.  
For the test data, a common strategy is to use the traces of special solutions such as CGO solutions.
And for certain cases without CGO solutions, other devices have been proposed. For example, 
``oscillating decaying solutions'' were developed to handle various types of anisotropic equations 
\cite{N2004_ODSapp, NUW2005_ODS, NUW2006, L2016,K2021,K2021-corrigendum, KLS2015};
"Wolff solutions" to the $p$-Laplace equation are employed in \cite{BKS2015,BHKS2018}.

However, instead of seeking a suitable class of solutions to \eqref{eq:main} directly, 
we tried to explore the possibility of using approximate 
solutions in this paper. This endeavor was motivated by the work \cite{LLLS2021} of Lassas, Liimatainen, Lin and Salo, where they 
exploit higher order linearizations of the Dirichlet-to-Neumann map to (among other things) solve the Calder\'{o}n type problems 
for \eqref{eq:main}. Following their idea, we will construct approximate solutions $\wt u_f$ for given Dirichlet data $f$ in such 
a way that the map $u\mapsto \wt u_f$ is the $m$-th Taylor polynomial of the 
(exact) solution operator. The primary innovation of this work is that we then
incorporate these approximate solutions to define an indicator functional, and
it is interesting that, by using this indicator functional, we can still use the 
Calder\'{o}n-type harmonic functions $\apa e^{\beta x\cdot(\oa+i\oa^\perp)}$ 
as test data to accomplish our reconstruction task. 

\subsection*{Organization}
The rest of the paper is organized as follows. In Section \ref{sec:prel} we give some 
preliminary results that will be needed, 
in particular a well-posedness result for \eqref{eq:main}. In Section \ref{sec:apps} 
we construct approximate solutions
obtained from Taylor approximation of the solution operator. In Section \ref{sec:indl}, 
we use these 
approximate solutions to define the indicator functional to be used, and give some
important observations. The precise definition of our test data 
and the main theorem is given in Section \ref{sec:recd}.

\section{Preliminaries}\label{sec:prel}

In this section, we give some preliminary results that will be used in 
subsequent sections. For simplicity let us introduce the following notation:
\begin{align*}
\ml A^2(\Oa) &=  H^2(\Oa)\cap L^\infty(\Oa),\\
\ml A^{3/2}(\pl\Oa) & = H^{3/2}(\pl\Oa)\cap L^\infty(\pl\Oa).
\end{align*}
We remark that all the function spaces considered are complex-valued.
$\ml A^2(\Oa)$ and $\ml A^{3/2}(\pl\Oa)$ are Banach spaces with the norms
\begin{align*}
\|\cdot\|_{\ml A^2(\Oa)}&:=\|\cdot\|_{H^2(\Oa)}+\|\cdot\|_{L^{\infty}(\Oa)},\\
\|\cdot\|_{\ml A^{3/2}(\pl\Oa)}&:=\|\cdot\|_{H^{3/2}(\pl\Oa)}
+\|\cdot\|_{L^{\infty}(\pl\Oa)}.
\end{align*}
Moreover, it is known (see Remark below) that 
$\ml A^2(\Oa)$ and $\ml A^{3/2}(\pl\Oa)$ are algebras with respect to multiplication. 
That is, $f,g\in\ml A^2(\Oa)$ (resp. $\ml A^{3/2}(\pl\Oa)$) implies $ fg\in\ml A^2(\Oa)$
(resp. $\ml A^{3/2}(\pl\Oa)$), and 
\begin{align}\label{algine}
\begin{aligned}
\|fg\|_{\ml A^2(\Oa)} &\le C \|f\|_{\ml A^2(\Oa)} \|g\|_{\ml A^2(\Oa)};\\
\|fg\|_{\ml A^{3/2}(\pl\Oa)} &\le C \|f\|_{\ml A^{3/2}(\pl\Oa)} 
\|g\|_{\ml A^{3/2}(\pl\Oa)},
\end{aligned}
\end{align}
where $C>0$ is independent of $f$ and $g$.

\begin{rem*}
Kato and Ponce \cite{KP1988} proved that $W^{s,p}(\mb R^n)\cap L^\infty(\mb R^n)$ 
is an algebra for $s>0$ (not necessarily an integer) and $1<p<\infty$. Based on this fact, that 
$H^{3/2}(\pl\Oa)\cap L^\infty(\pl\Oa)$ is an algebra 
is a straightforward consequence of the definition of $H^{3/2}(\pl\Oa)$ (see e.g. \cite{McLeanbook}) and the fact that
$L^\infty(\pl\Oa)$ is an algebra. On the other hand, 
as pointed out in \cite[Section 1.8.1]{Mazya2011},
$H^2(\Oa)\cap L^\infty(\Oa)$ is also an algebra since there exists 
a bounded linear operator from $H^2(\Oa)\cap L^\infty(\Oa)$ to 
$H^2(\mb R^n)\cap L^\infty(\mb R^n)$.
\end{rem*}

We are going to give a well-posedness result.
Motivated by \cite{HL2023}, we consider a more general form of nonlinear term $a(x,u)$, 
which is analytic in the second variable. Precisely, 
\begin{align}\label{eq:aform}
a(x,u)=\sum_{j= 2}^{\infty}\frac{a_{j}(x)}{j!}u^{j},\quad a_j\in L^\infty(\Oa),
\quad \sup_{j}\|a_j\|_{L^\infty(\Oa)}<\infty.
\end{align}
Note that there are no $a_0(x)$ and $a_1(x)$. 
The solution space to be considered will be the following subspace of $\ml A^2(\Oa)$:
\begin{align}\label{eq:def-V_s}
\ml V := \{ v\in\ml A^2(\Oa):\Da v\in L^{\infty}(\Oa)\mbox{ and }
v|_{\pl\Oa}\in L^{\infty}(\pl\Oa)\}.
\end{align}
$\ml V$ will be regarded as a Banach space with its own norm:
\begin{align*}
\| v\|_{\ml V}:=\|v\|_{\ml A^2(\Oa)}+\|\Da v\|_{L^{\infty}(\Oa)}+\|v|_{\pl\Oa}\|_{L^{\infty}(\pl\Oa)}.
\end{align*}

We now state the well-posedness result.
Note that similar results 
have been established in \cite{LLLS2021, HL2023, lu2022inverse, Nu2023}.

\begin{prop}\label{pro:well-posedness}
Let $a:\Oa\times\mb C\to\mb C$ satisfy \eqref{eq:aform}.
Then there exist $\da>0$ and $C>0$ such that for every 
\begin{align*}
f\in U_{\da}:=\{ f\in\ml A^{3/2}(\pl\Oa):\|f\|_{\ml A^{3/2}(\pl\Oa)}<\da\},
\end{align*}
there is a solution $u_f\in\ml V$ to the boundary value problem
\begin{align}\label{eq:main-a}
\begin{cases}
\Delta u+a(x,u)=0 & \mbox{ in }\Oa\\
u = f & \mbox{ on }\pl\Oa,
\end{cases}
\end{align}
which satisfies
\begin{align*}
\|u_{f}\|_{\ml V}\leq C\|f\|_{\ml A^{3/2}(\pl\Oa)}.
\end{align*}
Moreover, the solution operator 
$S:U_\da\to\ml V$, $f\mapsto u_f$, is a diffeomorphism between 
$U_\da$ and a 
neighborhood of $0$ in $\ml V$.
\end{prop}

Before giving the proof, recall that if a function $f:X\to Y$ between Banach spaces has $k$-th
(Fr\'{e}chet) derivative $(D^k f)_u$ at $u\in X$, then 
$(D^k f)_u$ can be regarded as a $k$-linear map from $X$ to $Y$ (see e.g. \cite[Section 1.1]{Hormander1}).
If $(D^k f)_v$ exists for $v$ in a neighborhood of $u$, then $(D^{k+1} f)_u$ is the unique $(k+1)$-linear
map $L$ from $X$ to $Y$ such that 
\begin{align*}
\sup\|((D^k f)_{u+h}-(D^k f)_u)(h_1,\ldots,h_k)- L(h,h_1,\ldots,h_k)\|_Y = o(\|h\|_X),
\end{align*}
where the $\sup$ is taken over all $h_1,h_2,\ldots,h_k\in X$ that satisfy $\|h_j\|_X=1$ for all $j$.

In the following lemma and its proof, we regard $a_0(x)=a_1(x)=0$ for convenience.
\begin{lem}\label{lem: a-diff} 
Let $a:\Oa\times \mb C\to \mb C$ satisfy \eqref{eq:aform}.
Then we have
$a(x,u(x))\in L^\infty(\Oa)$ for $u\in \ml V$. Moreover, 
\begin{align*}
G: \ml V \to L^\infty(\Oa),\quad u \mapsto a(x,u(x))
\end{align*}
is $C^{\infty}$ map, and its $k$-th derivative, $k\in\mb N$,
is given by
\begin{equation}\label{eq:G_ell-derivative}
(D^k G)_u(h_{1},h_{2},\cdots,h_{k})
=\left(\sum_{j= 0}^{\infty}\frac{a_{j+k}(x)}{j!}u^{j}\right)h_{1}h_{2}\cdots h_{k},
\end{equation}
for $u, h_1,\ldots,h_k\in \ml V$. 
\end{lem}

\begin{proof}
We will prove that $G: u\mapsto a(x,u(x))$ is indeed a $C^\infty$ map from $L^\infty(\Oa)$ to $L^\infty(\Oa)$, with derivatives given 
by the formula \eqref{eq:G_ell-derivative}, for $u,h_1,\ldots,h_k\in L^\infty(\Oa)$. Since 
$\ml V$ is continuously embedded in $L^\infty(\Oa)$, this implies the assertions of the lemma.

For simplicity we write $L^\infty$ instead of $L^\infty(\Oa)$ in the following.
Note that $L^\infty$ is an algebra, with $\|fg\|_{L^\infty}\le \|f\|_{L^\infty}\|g\|_{L^\infty}$ for $f,g\in L^\infty$.
By assumption, there is a positive constant $C_a$ independent of $j$ such that  
\begin{align*}
\|a_j\|_{L^\infty} \le C_a\ \,\mbox{for all}\ \, j.
\end{align*}
Thus, for all $u\in L^\infty$,
\begin{align*}
\|a(x,u(x))\|_{L^\infty}\le \sum_{j=0}^\infty \frac{C_a}{j!}\|u\|_{L^\infty}^j<\infty,
\end{align*}
and hence $G$ is a map from $L^\infty$ to $L^\infty$.

For $k\in\mb N$, define $G_k:L^\infty \to L^\infty$ by
\begin{align*}
G_k(u) = \sum_{j=0}^\infty \frac{a_{j+k}(x)}{j!}u^j.
\end{align*}
What we want to show is
\begin{align}\label{eq:tsla}
(D^k G)_u(h_1,\ldots,h_k) = G_k(u)h_1\cdots h_k\quad(u,h_1,\ldots,h_k \in L^\infty).
\end{align}
For $u,h\in L^\infty$, we have
\begin{align*}
G(u+h)-G(u)
&=  \sum_{j=1}^{\infty}\frac{a_{j}(x)}{j!}\left((u+h)^{j}-u^{j}\right)\\
&=  \sum_{j=1}^{\infty}\frac{a_{j}(x)}{j!}\sum_{\ell=1}^{j}\begin{pmatrix}j\\
\ell
\end{pmatrix}u^{j-\ell}h^{\ell}\\
&= \sum_{\ell=1}^\infty\left(\sum_{j=\ell}^\infty \frac{a_{j}(x)}{(j-\ell)!}u^{j-\ell}\right)\frac{h^\ell}{\ell!}\\
&= \sum_{\ell=1}^\infty G_\ell(u)\frac{h^\ell}{\ell!}.
\end{align*}
Thus we obtain
\begin{align*}
G(u+h)-G(u)-G_1(u)h = \sum_{\ell=2}^\infty G_\ell(u)\frac{h^\ell}{\ell !}.
\end{align*}
It is easy to see that $\|G_k(u)\|_{L^\infty} \le C_a e^{\|u\|_{L^\infty}}$ for all $k\in\mb N$, and hence the above formula implies 
\begin{align*}
\|G(u+h)-G(u)-G_1(u)h \|_{L^\infty} = O(\|h\|_{L^\infty}^2)\quad\mbox{as}\quad \|h\|_{L^\infty}\to 0.
\end{align*}
This proves \eqref{eq:tsla} for $k=1$. Note that we can regard $G$ as $G_0$, and the above equation says 
$(DG_0)_u(h)=G_1(u)h$.
In the same manner, we can prove that 
\begin{align}\label{eq:obv}
(D G_k)_u(h) = G_{k+1}(u)h\quad (u,h\in L^\infty)
\end{align}
for all $k\in\mb N$.

To establish \eqref{eq:tsla} for $k \ge 2$, we proceed by induction. Suppose we have shown that it is 
true for some $k\in\mb N$. Then, for $u,h\in L^\infty$ and $h_1,\ldots,h_k\in L^\infty$ such that $\|h_j\|_{L^\infty}=1$ for all $j$, 
we have
\begin{align*}
\Big\|\Big((D^k &G)_{u+h} - (D^k G)_u\Big)(h_1,\ldots,h_k) - G_{k+1}(u)hh_1\ldots h_k\Big\|_{L^\infty}\\
& = \Big\|\Big(G_{k}(u + h) - G_{k}(u)\Big)h_1\cdots h_k - G_{k+1}(u)hh_1\ldots h_k\Big\|_{L^\infty} \\
& \le  \Big\|\Big(G_{k}(u + h) - G_{k}(u)\Big) - G_{k+1}(u)h\Big\|_{L^\infty}.
\end{align*}
From \eqref{eq:obv}, the last expression is $o(\|h\|_{L^\infty})$, and the proof is completed.
\end{proof}

We need one more lemma.
\begin{lem}\label{lem:elliptic_reg}
For any $(g,f)\in L^\infty(\Oa)\times\ml A^{3/2}(\pl\Oa)$, the boundary value problem
\begin{equation}\label{lee}
\begin{cases}
\Delta v=g & \mbox{ in }\Oa\\
v=f & \mbox{ on }\pl\Oa
\end{cases}
\end{equation}
has a unique solution $v\in\ml V$, which satisfies
\begin{align*}
\|v\|_{\ml V}\le C( \|g\|_{L^\infty(\Oa)} + \|f\|_{\ml A^{3/2}(\pl\Oa)} ),
\end{align*}
where $C>0$ is a constant independent of $g$ and $f$.
\end{lem}
\begin{proof}
It is well-known that \eqref{lee} has a unique solution $v\in H^2(\Oa)$, 
which satisfies 
\begin{align*}
\|v\|_{H^2(\Oa)} 
&\le C (\|g\|_{L^{2}(\Oa)}+\|f\|_{H^{3/2}(\pl\Oa)})\\
&\le C (\|g\|_{L^{\infty}(\Oa)}+\|f\|_{H^{3/2}(\pl\Oa)}).
\end{align*}
On the other hand, by \cite[Lemma 1.14]{FR2022}, we also have
\begin{align*}
\|v\|_{L^{\infty}(\Oa)} \leq C (\|g\|_{L^{\infty}(\Oa)}+\|f\|_{L^{\infty}(\pl\Oa)}).
\end{align*}
Thus, by definition,
\begin{align*}
\|v\|_{\ml V} 
&= \|v\|_{H^2(\Oa)} + \|v\|_{L^\infty(\Oa)} + \|g\|_{L^\infty(\Oa)} + \|f\|_{L^\infty(\pl\Oa)}\\
&\le C\big(\|g\|_{L^\infty(\Oa)}+ \|f\|_{\ml A^{3/2}(\pl\Oa)}\big).
\end{align*}
\end{proof}

We now give the proof of Proposition \ref{pro:well-posedness}. The idea of using implicit function theorem follows \cite{LLLS2021}.

\begin{proof}[Proof of Proposition \ref{pro:well-posedness}]
Consider the map
\begin{gather*}
F:\ml A^{3/2}(\pl\Oa)\times \ml V
\to L^\infty(\Oa)\times\ml A^{3/2}(\pl\Oa),\\
F(f,u)=(\Da u+a(x,u),\,u|_{\pl\Oa}-f).
\end{gather*}
Clearly, $F(0,0)=(0,0)$. Notice that except for $a(x,u)$, the other parts of $F$ are linear. 
Hence, by Lemma \ref{lem: a-diff}, $F$ is a $C^{\infty}$ map, and the partial derivative of $F$ at 
$(0,0)$ in the $u$-variable is given by
\[
(D_{u}F)_{(0,0)}(v)=(\Da v+(DG)_{0}(v),v|_{\pl\Oa})=(\Da v,v|_{\pl\Oa}).
\]
By Lemma \ref{lem:elliptic_reg}, $(D_{u}F)_{(0,0)}$ is a Banach space isomorphism from $\ml V$ to 
$L^\infty(\Oa)\times\ml A^{3/2}(\pl\Oa)$. Therefore, the implicit function theorem for Banach spaces (see e.g. \cite[Theorem 10.6]{Renardybook}) implies that 
there exists a $\da$-neighborhood $U_\da$ of $0$ in $\ml A^{3/2}(\pl\Oa)$, i.e.
\begin{align*}
U_{\da}=\big\{ f\in\ml A^{3/2}(\pl\Oa):
\|f\|_{\ml A^{3/2}(\pl\Oa)}<\da\big\},
\end{align*}
and a $C^1$ map $S:U_{\da}\to \ml V$ such that $F(f,S(f))=0$ for all 
$f\in U_\da$. 
Moreover, since $F$ is smooth, by further restricting the size of $\da$ if necessary, the following can be 
guaranteed: 
\begin{enumerate}
\item $S$ is a diffeomorphism between $U_\da$ and a neighborhood of $0$ in 
$\ml V$.
\item $S$ is Lipschitz continuous on $U_\da$. Hence, for all $f\in U_\da$, 
\begin{align*}
\|S(f)\|_{\ml V} = \|S(f)-S(0)\|_{\ml V}\le C\|f-0\|_{\ml A^{3/2}(\pl\Oa)} 
= C\|f\|_{\ml A^{3/2}(\pl\Oa)}
\end{align*}
for some constant $C>0$ independent of $f$.
\end{enumerate}
Thus, by letting $u_f= S(f)$, $u_f$ satisfies the assertions of the 
proposition.
\end{proof}

\section{The Indicator Functional}

In this section, we begin our discussion on the application of the enclosure method to \eqref{eq:main}.
As mentioned in the introduction, we will construct approximate solutions to 
\eqref{eq:main}, and then use them to define an indicator functional. 

\subsection{Taylor approximation of the solution operator}\label{sec:apps}

Recall our equation from Section \ref{sec:prob}.
By Proposition \ref{pro:well-posedness}, there exist $\da>0$ and $S\in C^\infty(U_\da,\ml V)$ 
such that for $f\in U_\da$, $u_f=S(f)$ is the unique solution of
\begin{align}\label{eq:mbvp}
\begin{cases}
\Da u + q(x)u^m = 0 & \mbox{in }\Oa\\
u = f & \mbox{on }\pl\Oa
\end{cases}
\end{align}
in a neighborhood of $0\in\ml V$. We now follow the idea of \cite{LLLS2021} to derive the Taylor approximation for $S$.
In fact, since $S$ is smooth and $S(0)=0$, we expect that $S(f)$
can be approximated by Taylor's polynomial:
\begin{align}\label{otr}
S(f)\approx\sum_{j=0}^{k}\frac{(D^{j}S)_{0}(f,\cdots,f)}{j!}\qquad (k\in\mb N).
\end{align}
The key observation here is that the monomials in \eqref{otr} can be obtained by differentiating the curve
$\vn \mapsto S(\vn f)=u_{\vn f}$ ($-1<\vn<1$) at $\vn =0$. Precisely, let $d_\vn^j$ denotes $(d/d\vn)^j$. We have
\begin{align*}
(D^j S)_0(f,\cdots,f) = d_\vn^j u_{\vn f}|_{\vn=0},
\end{align*}
which will temporarily be denoted by $u^{(j)}_f$. Thus \eqref{otr} can be rewritten as 
\begin{align}\label{ntr}
S(f)\approx\sum_{j=0}^k u^{(j)}_f/j!.
\end{align}
By Taylor's theorem, one expects that the approximation \eqref{ntr} has a remainder which is 
$O(\|f\|^{k+1}_{\ml A^{3/2}(\pl\Oa)})$. Nevertheless, we are going to show that, for $k=m$, 
\eqref{ntr} really provides a better approximation than expected, with an 
$O(\|f\|_{\ml A^{3/2}(\pl\Oa)}^{2m-1})$ remainder.

In the following we (formally) deduce the boundary value problems satisfied by 
$u^{(j)}_f$, $j=1,\ldots,m$. (Note that $u_f^{(0)} = u_{0f}=S(0)=0$.)
By definition, $u_{\vn f}$ satisfies 
\begin{align}\label{eq:0th}
\begin{cases}
\Da u_{\vn f} + q(x)u_{\vn f}^m = 0 & \mbox{in }\Oa\\
u_{\vn f} = \vn f & \mbox{on }\pl\Oa.
\end{cases}
\end{align}
Applying $d_\vn^j(\cdot)|_{\vn =0}$ to \eqref{eq:0th}, we get
\begin{align}\label{eq:1stp}
\begin{cases}
\Da u_f^{(j)} + q(x)d_\vn^j(u_{\vn f}^m)|_{\vn=0} = 0 & \mbox{in }\Oa\\
u_f^{(j)} = d_\vn^j (\vn f)|_{\vn=0} & \mbox{on }\pl\Oa.
\end{cases}
\end{align}
Since  $u_{0 f} = S(0) =0$, for $j=1$ the differentiation of the nonlinear term is 
\begin{align*}
d_\vn(u_{\vn f}^m)|_{\vn=0} = m u_{0 f}^{m-1} u_f^{(1)}=0.
\end{align*}
And the differentiation of the boundary condition is
\begin{align*}
d_\vn^j (\vn f)|_{\vn=0} = f|_{\vn=0} = f.
\end{align*}
Thus we obtain
\begin{align}\label{eq:1st}
\begin{cases}
\Da u_f^{(1)} = 0 & \mbox{in }\Oa\\
u_f^{(1)} = f & \mbox{on }\pl\Oa.
\end{cases}
\end{align} 
In fact, by the same reason that $u_{0 f}=0$, we have
\begin{align*}
d_\vn^j (u_{\vn f}^m)|_{\vn = 0} = 0\quad\mbox{for}\quad 2\le j\le m-1,
\end{align*}
and
\begin{align*}
d_\vn^m (u_{\vn f}^m)|_{\vn = 0} = m! (u^{(1)}_f)^m.
\end{align*}
On the other hand, $d_\vn^j(\vn f) = 0$ for all $j\ge 2$. Hence, for $2\le j\le m-1$, 
\eqref{eq:1stp} gives
\begin{align*}
\begin{cases}
\Da u_f^{(j)} = 0 & \mbox{in }\Oa\\
u_f^{(j)} = 0 & \mbox{on }\pl\Oa,
\end{cases}
\end{align*}
which implies $u_f^{(j)}=0$. And for $j=m$, \eqref{eq:1stp} gives
\begin{align}\label{eq:mth}
\begin{cases}
\Da u_f^{(m)} = -m! q(x) (u^{(1)}_f)^m & \mbox{in }\Oa\\
u_f^{(m)} = 0 & \mbox{on }\pl\Oa.
\end{cases}
\end{align}

In conclusion, the right-hand side of \eqref{ntr} with $k=m$ is given by 
\begin{align*}
\wt u_f := u_f^{(1)}+u_f^{(m)}/m!.
\end{align*}
By \eqref{eq:1st} and \eqref{eq:mth}, $\wt u_f$ satisfies 
\begin{align*}
\begin{cases}
\Da \wt u_f = -q(x) (u_f^{(1)})^m & \mbox{in }\Oa,\\
\wt u_f = f & \mbox{on }\pl\Oa.
\end{cases}
\end{align*}
To avoid cumbersome notation, we will write $v_f$ for $u_f^{(1)}$.
By Lemma \ref{lem:elliptic_reg}, we can redefine $v_f$ and $\wt u_f$ as follows.

\begin{defn}\label{defvtu}
For $f\in U_\da$, let $v_f$ be the unique function in $\ml V$ that satisfies 
\begin{align*}
\begin{cases}
\Da v_f = 0 &\mbox{in }\Oa,\\
v_f = f &\mbox{on }\pl\Oa.
\end{cases}
\end{align*}
In other words, $v_f$ is the unique harmonic function in $\Oa$ with trace $f$.
And then let $\wt u_f$ be the unique function in $\ml V$ that satisfies 
\begin{align*}
\begin{cases}
\Da \wt u_f = - q(x)v_f^m &\mbox{in }\Oa,\\
\wt u_f = f &\mbox{on }\pl\Oa.
\end{cases}
\end{align*}
\end{defn}
 
We now give the error estimates.

\begin{lem}\label{lem:est12}
For $f\in U_\da$, we have
\begin{align}
&\|u_f - v_f\|_{\ml V} \le C \|f\|_{\ml A^{3/2}(\pl\Oa)}^m, \label{est1}\\
&\|u_f - \wt u_f \|_{\ml V} \le C \|f\|_{\ml A^{3/2}(\pl\Oa)}^{2m-1}, \label{est2}
\end{align}
where $C>0$ is a constant independent of $f$.
\end{lem}

\begin{proof}
By \eqref{eq:mbvp} and Definition \ref{defvtu}, 
\begin{align*}
\begin{cases}
\Da (u_f - v_f) = -q(x)u_f^m & \mbox{in }\Oa\\
u_f -  v_f = 0 & \mbox{on }\pl\Oa.
\end{cases}
\end{align*}
By Lemma \ref{lem:elliptic_reg} and Proposition \ref{pro:well-posedness},
\begin{align*}
\|u_f - v_f\|_{\ml V} \le C\|q(x)u_f^m\|_{L^\infty(\Oa)} \le C \|u_f\|_{L^\infty(\Oa)}^m
\le C\|f\|_{\ml A^{3/2}(\pl\Oa)}^m,
\end{align*}
and \eqref{est1} is proved. Similarly, we have
\begin{align*}
\begin{cases}
\Da (u_f - \wt u_f) = -q(x)\big(u_f^m - v_f^m\big) & \mbox{in }\Oa\\
u_f - \wt u_f= 0 & \mbox{on }\pl\Oa,
\end{cases}
\end{align*}
and Lemma \ref{lem:elliptic_reg} implies
\begin{align*}
\|u_f - \wt u_f\|_{\ml V} 
&\le C\|q(x)(u_f^m - v_f^m)\|_{L^\infty(\Oa)} \\
&\le C\|u_f- v_f\|_{L^\infty(\Oa)}
\sum_{j=0}^{m-1} \|u_{f}\|_{L^\infty(\Oa)}^{m-1-j}\|v_{f}\|_{L^\infty(\Oa)}^{j}.
\end{align*}
By the just proved inequality \eqref{est1}, $\|u_f- v_f\|_{L^\infty(\Oa)}\le C\|f\|_{\ml A^{3/2}(\pl\Oa)}^m$;
By Proposition \ref{pro:well-posedness}, $\|u_f\|_{L^\infty(\Oa)}\le \|f\|_{\ml A^{3/2}(\pl\Oa)}$; and 
by applying Lemma \ref{lem:elliptic_reg} to \eqref{eq:1st}, we also have 
$\|v_f\|_{L^\infty(\Oa)}\le \|f\|_{\ml A^{3/2}(\pl\Oa)}$. Thus the above inequality implies \eqref{est2}.
\end{proof}

We will also need the analogue of $\wt u_{f}$ corresponding to the background coefficient $q_0(x)$, which 
we will denote by $\wt u_{0,f}$. Precisely, we make the following definition. 
(By contrast, note that $v_f$ is independent of $q(x)$, and there is no need to consider $v_{0,f}$.)

\begin{defn}\label{ut0f}
Let $v_f$ be as defined in Definition \ref{defvtu}. Then, let $\wt u_{0,f}$ be the unique function in 
$\ml V$ that satisfies 
\begin{align*}
\begin{cases}
\Da \wt u_{0,f} = - q_0(x)v_f^m &\mbox{in }\Oa,\\
\wt u_{0,f} = f &\mbox{on }\pl\Oa.
\end{cases}
\end{align*}
\end{defn}

\subsection{Definition of the indicator functional}\label{sec:indl}

In the following let us write $\pl_\nu u$ for $(\pl_\nu u)|_{\pl\Oa}$.
Note that for $u\in H^2(\Oa)$ it is a formal expression of $\mbox{tr}(\na u)\cdot \nu$,
where $\mbox{tr}: H^{1}(\Oa)\to H^{1/2}(\pl\Oa)$ 
is the trace operator and $\nu$ is the outward normal on $\pl\Oa$. 
It is easy to check that 
\begin{align}\label{sti}
\|\pl_\nu u\|_{H^{1/2}(\pl\Oa)}\le C\|u\|_{H^2(\Oa)},
\end{align}
where $C>0$ is independent of $u$.

The indicator functional we are going to use is
\begin{equation}\label{eq:E-def}
E(f) := \int_{\pl\Oa}(\pl_{\nu}u_{f}-\pl_{\nu}\wt u_{0,f})\ol{f^{m}}\,d\sigma \qquad(f\in U_\da).
\end{equation}
In addition to $E(f)$, we also define the ``auxiliary indicator functional'':
\begin{align}\label{eq:et-def}
\wt E(f):=\int_{\pl\Oa}(\pl_{\nu}\wt u_{f}-\pl_{\nu}\wt u_{0,f})\ol{f^{m}}\,d\sigma\qquad(f\in U_\da).
\end{align}
Note that since $q$ is unknown, we are not able to obtain $\pl_{\nu}\wt u_{f}$,
and hence 
using $\wt E(f)$ to reconstruct $D$ is not realistic. Nevertheless, the following lemma shows 
that $\wt E(f)$ is very close to $E(f)$ for small $f$. It will 
allow us to study $E(f)$ by analyzing $\wt E(f)$.

\begin{lem}\label{lem:approx_indicator}
For $f\in U_{\da}$, we have
\[
|E(f)-\wt E(f)|\leq C\|f\|_{\ml A^{3/2}(\pl\Oa)}^{3m-1},
\]
where $C>0$ is independent of $f$.
\end{lem}
\begin{proof}
We have
\begin{align*}
|E(f)-\wt E(f)| & \le \|\pl_{\nu}u_{f}-\pl_{\nu}\wt u_{f}\|_{L^2(\pl\Oa )}\|f^m\|_{L^2(\pl\Oa)}\\
 & \leq\|\pl_{\nu}u_{f}-\pl_{\nu}\wt u_{f}\|_{H^{1/2}(\pl\Oa)}\|f^{m}\|_{\ml A^{3/2}(\pl\Oa)}\\
 & \leq C\|u_{f}-\wt u_{f}\|_{H^{2}(\Oa)}\|f\|_{\ml A^{3/2}(\pl\Oa)}^m\\
 & \leq C\|f\|_{\ml A^{3/2}(\pl\Oa)}^{3m-1}.
\end{align*}
The third inequality is true by \eqref{sti} and \eqref{algine}, and the last inequality is given
by Lemma \ref{lem:est12}.
\end{proof}

In the following we apply Green's identity to rewrite $\wt E(f)$ ($f\in U_\da$) in a different form,
and then give an important observation (which 
actually is the reason for defining our indicator functional in the above form.)
For this, first recall that $q(x)=q_0(x)+\chi_D q_D(x)$. Thus, 
from Definition \ref{defvtu} and Definition \ref{ut0f}, $\wt u_f - \wt u_{0,f}$ satisfies 
\begin{align}\label{wtz}
\begin{cases}
\Da (\wt u_f- \wt u_{0,f}) = -\chi_D q_D(x) v_f^m & \mbox{in }\Oa\\
\wt u_f - \wt u_{0,f} = 0 & \mbox{on }\pl\Oa.
\end{cases}
\end{align}
On the other hand, note that $f^m\in\ml A^{3/2}(\pl\Oa)$ since $\ml A^{3/2}(\pl\Oa)$ is an algebra.
Hence $f^m = F_{m,f}|_{\pl\Oa}$ for some function $F_{m,f}\in H^2(\Oa)$.
By Green's second identity, we then obtain
\begin{align}
\wt E(f) 
& = \int_{\Oa}\ol{F_{m,f}}\Da(\wt u_f - \wt u_{0,f})\,dx 
- \int_{\Oa}(\wt u_f - \wt u_{0,f}) \ol{\Da F_{m,f}} \,dx \nonumber\\
&\quad +\int_{\pl\Oa} (\wt u_f - \wt u_{0,f})\ol{\pl_\nu F_{m,f}}\,d\sigma \nonumber\\
& = -\int_{D}q_D(x) v_f^m \ol{F_{m,f}} \,dx - \int_\Oa (\wt u_f - \wt u_{0,f})\ol{\Da F_{m,f}} \,dx. \label{prest}
\end{align}
At this point, we would like to choose $F_{m,f}=v_f^m$ so that the first integral 
in \eqref{prest} will provide a 
``lower bound estimate'' for $|\wt E(f)|$, which is crucial in proving the reconstruction theorem.
A possible problem is that in general we do not know if $v_f^m$ satisfies 
\begin{align}\label{sbtri}
v_f^m|_{\pl\Oa} = (v_f|_{\pl\Oa})^m = f^m.
\end{align}
Notice however that if $v_f\in C(\ol{\Oa})$ and hence 
$v_f|_{\pl\Oa}=f$ in the classical sense, then \eqref{sbtri} holds true obviously. 
On the other hand, for the second integral, it would be preferable if we could eliminate it.
Fortunately, by considering $f$ to be the traces 
of Calder\'{o}n-type harmonic functions, all the 
desires can be fulfilled. Precisely, we have the following result. 
\begin{lem}\label{lem:wtee}  
Let $f(x)=a e^{b x\cdot(\oa+i\oa^\perp)}|_{\pl\Oa}$,
where $a,b$ are constants, and $\oa,\oa^\perp\in\mb S^{n-1}$ 
are such that $\oa\perp\oa^\perp$. Then $v_f(x) = a e^{b x\cdot(\oa+i\oa^\perp)}$ (for $x\in\Oa$), and 
\begin{align}\label{thest}
\wt E(f) = -\int_D q_D(x) |v_f|^{2m}\,dx.
\end{align}
\end{lem}
\begin{proof}
Since $a e^{b x\cdot(\oa+i\oa^\perp)}$ is a harmonic function whose trace on $\pl\Oa$ equals
$f$, it must be $v_f$ (in $\Oa$). Then, observe that
$v_f^m =a^m e^{mb x\cdot(\oa+i\oa^\perp)}$ is still a harmonic function. In fact, 
it is the unique harmonic function whose trace on $\pl\Oa$  is $f^m$.
As a consequence, we can set $F_{m,f} = v_f^m$ in \eqref{prest}, which gives \eqref{thest}.
\end{proof}

\section{Reconstruction of $D$}\label{sec:recd}

We give our reconstruction theorem in this section. 
First, set the test data as follows: Given $h>0$, $J>0$, $t\in\mb R$, and $\oa,\oa^\perp\in \mb S^{n-1}$ such
that $\oa\perp\oa^\perp$, let 
\begin{align}\label{dovf}
\begin{aligned}
v_h(x) &:= e^{-\frac{J}{h}}e^{-\frac{1}{h}(x\cdot\oa+ix\cdot\oa^{\perp}-t)}|_\Oa,\\
f_h(x) &:= e^{-\frac{J}{h}}e^{-\frac{1}{h}(x\cdot\oa+ix\cdot\oa^{\perp}-t)}|_{\pl\Oa}.
\end{aligned}
\end{align}
Note that $v_h= v_{f_h}$. Let
\begin{align}\label{ians}
b = b(\oa) := \inf_{x\in\Oa}x\cdot \oa,\quad B = B(\oa):=\sup_{x\in\Oa}x\cdot \oa.
\end{align}
We will investigate whether the half-space $\{x\cdot\oa\le t\}$ intersects $D$ or not by
analyzing the limiting behavior of $E(f_h)$ as $h\to 0$. 
For this, it suffices to consider $b<t<B$.
On the other hand, remember that $E(f_h)$ is defined 
only for $f_h\in U_\da$. This requirement can be fulfilled by setting $J$ large enough.
In fact, we have the following lemma.

\begin{lem}\label{lem:calculate-f}
$\|f_h\|_{\ml A^{3/2}(\pl\Oa)} = O(h^{-2}e^{-\fc{J}{h}}e^{\frac{t-b}{h}})$
as $h\to 0$.
As a consequence, if $J>B-b$ and $t<B$, we have $f_h\in U_\da$ for all $h$ small enough. 
\end{lem}
\begin{proof}
For an $n$-dimensional multi-index $\ba=(\ba_1,\ldots,\ba_n)$, we have 
\begin{align*}
\pl^{\ba}v_{h}(x)
=(-1)^{|\ba|}h^{-|\ba|}e^{-\fc{J}{h}}
e^{-\frac{1}{h}(x\cdot\oa+ix\cdot\oa^{\perp}-t)}(\oa+i\oa^\perp)^{\ba},
\end{align*}
where for a vector $a=(a_1,\ldots,a_n)$, $a^\ba:=\prod_{j} a_j^{\ba_j}$.
It's not hard to check that $|(\oa+i\oa^\perp)^{\ba}|\le 1$.
Thus 
\begin{align*}
\|\pl^{\ba }v_h\|_{L^\infty(\Oa)}
\le  h^{-|\ba |}e^{-\fc{J}{h}}e^{\frac{t-b}{h}}.
\end{align*}
Since $\Oa$ is bounded, 
\begin{align}\label{l2est}
\|\pl^{\ba }v_h\|_{L^2(\Oa)}\le |\Oa|^{1/2}\|\pl^{\ba}v_h\|_{L^\infty(\Oa)}\le |\Oa|^{1/2}h^{-|\ba |}e^{-\fc{J}{h}}e^{\frac{t-b}{h}}.
\end{align}
By \eqref{l2est} and the boundedness of the 
trace operator from $H^2(\Oa)$ to $H^{3/2}(\pl\Oa)$, we obtain
\begin{align*}
\|f_h\|_{H^{3/2}(\pl\Oa)} \le C\|v_h\|_{H^2(\Oa)} = O(h^{-2}e^{-\fc{J}{h}}e^{\frac{t-b}{h}})
\quad\mbox{as}\quad h\to 0.
\end{align*}
On the other hand, by definition of $f_h$ we have 
\begin{align*}
\|f_h\|_{L^\infty(\pl\Oa)}\le e^{-\frac{J}{h}}e^{\frac{t-b}{h}}.
\end{align*}
Hence $\|f_h\|_{\ml A^{3/2}(\pl\Oa)} = O(h^{-2}e^{-\fc{J}{h}}e^{\frac{t-b}{h}})$ as $h\to 0$.
\end{proof}

We are now ready to give our main theorem. Recall our equation and assumptions 
in Section \ref{sec:prob}, 
the definition of $E$ in \eqref{eq:E-def}, and the definition of $b$ and $B$ in \eqref{ians}.
Moreover, define 
\begin{align*}
t_* = \inf_{x\in D}x\cdot \oa.
\end{align*}

\begin{thm}\label{thm:main_thm} Let $f_h$ be defined by \eqref{dovf}, with 
$J>\frac{3m-1}{m-1}(B-b)$, $b<t<B$, and $\oa,\oa^\perp\in \mb S^{n-1}$ being such that $\oa\perp\oa^\perp$. Then we have the following assertions:
\begin{enumerate}
\item[\textup{(i)}]  If $\{ x\cdot\oa\leq t\}\cap \ol D =\emptyset$, that is, $t<t_*$, then 
\[
e^{\frac{2mJ}{h}}|E(f_{h})|\to0,\mbox{ as }h\to 0.
\]
\item[\textup{(ii)}]  If $\{ x\cdot\oa\le t\}\cap D \neq\emptyset$, that is, $t>t_*$, then
\[
e^{\frac{2mJ}{h}}|E(f_{h})|\to\infty,\mbox{ as }h\to 0.
\]
\end{enumerate}
\end{thm}

\begin{rem*}
In the above theorem, we do not consider the case that the half-space 
$\{x\cdot \oa \le t\}$ ``just touches'' the boundary 
of $D$, that is, $t = t_*$.
The analysis of this case is more complicated and would require a higher regularity assumption on $\pl D$. 
(By contrast, we only assume $D$ is an open set in the theorem.)
From the point of view of application, the above theorem already looks satisfactory. 
\end{rem*}

\begin{proof}[Proof of Theorem \ref{thm:main_thm}]
In this proof we assume that $h$ is small enough so that $f_h\in U_\da$. 
By Lemma \ref{lem:approx_indicator}, we have 
\begin{align}\label{ebyet}
\left|e^{\frac{2mJ}{h}}|E(f_h)|-e^{\frac{2mJ}{h}}|\wt E(f_h)|\right|
\le Ce^{\frac{2mJ}{h}}\|f_h\|_{\ml A^{3/2}(\pl\Oa)}^{3m-1}.
\end{align}
By Lemma \ref{lem:calculate-f}, $\|f_h\|_{\ml A^{3/2}(\pl\Oa)} \le C h^{-2}e^{-\fc{J}{h}}e^{\frac{t-b}{h}}$ for $h$ small, and hence
\begin{align*}
e^{\frac{2mJ}{h}}\|f\|_{\ml A^{3/2}(\pl\Oa)}^{3m-1} & \leq Ce^{\frac{2mJ}{h}}h^{-2(3m-1)}e^{-\frac{J}{h}(3m-1)}e^{\frac{t-b}{h}(3m-1)}\\
 & \leq Ce^{\frac{-J(m-1)+(t-b)(3m-1)}{h}}h^{-2(3m-1)}.
\end{align*}
Since $J>\frac{3m-1}{m-1}(B-b)$ and $b<t<B$, 
\begin{equation*}
e^{\frac{2mJ}{h}}\|f\|_{\ml A^{3/2}(\pl\Oa)}^{3m-1}\to 0 \mbox{ as }h \to 0.
\end{equation*}
Therefore, \eqref{ebyet} implies
\begin{align*}
\lim_{h\to 0}e^{\frac{2mJ}{h}}|E(f_{h})| = \lim_{h\to 0}e^{\frac{2mJ}{h}}|\wt E(f_{h})|,
\end{align*} 
including the possibility of ``$\infty=\infty$''. 
It remains to study the limit of $e^{\frac{2mJ}{h}}|\wt E(f_{h})|$ for the two cases considered 
in the theorem.

(i) $\{ x\cdot\oa\leq t\}\cap \ol D =\emptyset$, that is, $t<t_*$. 
Since $v_h = v_{f_h}$, by Lemma \ref{lem:wtee} we have the following upper-bound estimate:
\begin{align*}
e^{\frac{2mJ}{h}}|\wt E(f_{h})|
&=  e^{\frac{2mJ}{h}}\left|\int_{D}q_{D}|v_{h}|^{2m}\,dx\right|\\
&\leq  e^{\frac{2mJ}{h}}\|q_D\|_{L^\infty(\Oa)}\int_{D}e^{-\frac{2mJ}{h}}e^{-\frac{2m}{h}(x\cdot\oa-t)}\,dx \\
&\leq  \|q_D\|_{L^\infty(\Oa)}\int_{D}e^{-\frac{2m}{h}(t_{*}-t)}\,dx\to 0 \mbox{ as }h \to 0.
\end{align*}
Hence the assertion is true. 

(ii) $\{ x\cdot\oa\le t\}\cap D \neq\emptyset$, that is, $t>t_*$. In this case, we can choose 
a small enough $\vn>0$ such that $t-\vn$ is still greater than $t_*$. Thus 
$\{ x\cdot\oa\le t-\vn \}\cap D \neq\emptyset$. Since $D$ is open, this intersection obviously has
a positive measure, which we denote by $|\{ x\cdot\oa\le t-\vn \}\cap D|$.
By the jump condition
\eqref{eq:jump}, we have the following lower bound estimate:
\begin{align*}
e^{\frac{2mJ}{h}} |\wt E(f_{h})|
&=  e^{\frac{2mJ}{h}}\left|\int_{D}q_{D}|v_{h}|^{2m}\,dx\right|\\
&\geq  e^{\frac{2mJ}{h}}\mu \int_{D}e^{-\frac{2mJ}{h}}e^{-\frac{2m}{h}(x\cdot\oa-t)}\,dx\\
&\geq  \mu \int_{\{x\cdot\oa \le t-\vn \}\cap D}e^{-\frac{2m}{h}(x\cdot\oa-t)}\,dx\\
&\geq  \mu |\{x\cdot\oa \le t-\vn\}\cap D| e^{\frac{2m\vn}{h}}\to \infty \mbox{ as }h\to 0.
\end{align*}
Hence the assertion is also true. 
\end{proof}

\section*{Acknowledgments}
This work is partly supported by the National Science and Technology Council, Taiwan, under project 
number MOST 111-2115-M-006 -015 -MY2.

\bibliographystyle{plain}
\bibliography{ref}
 
\end{document}